\documentclass{amsart}
\usepackage{amsmath,amssymb,graphics}
\def\sumprime_#1{\setbox0=\hbox{$\scriptstyle{#1}$}
\setbox2=\hbox{$\displaystyle{\sum}$}
\setbox4=\hbox{${}'
\mathsurround=0pt$}
\dimen0=.5\wd0 \advance\dimen0 by-.5\wd2
\ifdim\dimen0>0pt
\ifdim\dimen0>\wd4 \kern\wd4 \else\kern\dimen0\fi\fi
\mathop{{\sum}'}_{\kern-\wd4 #1}}
\def\sumflat_#1{\setbox0=\hbox{$\scriptstyle{#1}$}
\setbox2=\hbox{$\displaystyle{\sum}$}
\setbox4=\hbox{${}\flat
\mathsurround=0pt$}
\dimen0=.5\wd0 \advance\dimen0 by-.5\wd2
\ifdim\dimen0>0pt
\ifdim\dimen0>\wd4 \kern\wd4 \else\kern\dimen0\fi\fi
\mathop{{\sum}^\flat}_{\kern-\wd4 #1}}

\begin{document}
\title{Are there infinitely many twin primes ?}
\author{D. A. Goldston}
\address{Department of Mathematics, San Jose
State University, San Jose, CA 95192, USA}
\email{goldston@math.sjsu.edu}

\maketitle

\section{Introduction}

Most people reading this article will recognize the sequence 
\begin{equation} 2,3,5,7,11,13,17,19, 23, 29,31,37,41,43,47,53, 61,67, 71,73,79,83,89,97 \ldots \label{1}\end{equation}
of prime numbers. The primes are defined to be the natural or counting numbers with exactly two factors, namely 1 and themselves, or equivalently the numbers only divisible by themselves and 1.\footnote{ While one could argue that 1 itself should be a prime, most mathematicians prefer to put 1 in a class by itself.} Here are the first 10 primes starting at a billion:
\begin{equation}  \begin{split}&1000000007, 1000000009, 1000000021, 1000000033, 1000000087, \\&1000000093, 
1000000097, 1000000103, 1000000123, 1000000181 .  \end{split} \label{2}  \end{equation}
You can easily find these primes yourself on your own computer using a mathematical package such as Mathematica or Maple. Further, no matter where you start looking, even at numbers with hundreds of digits,  you will have no trouble finding primes. Based on this experimental evidence, we can safely and scientifically conclude that the sequence of primes never ends.

\bigskip 
\noindent\textbf{Scientific or Empirical Observation.}  \emph{There are infinitely many prime numbers.} 
\bigskip

Is this observation a true fact? Certainly it is easy to verify that no matter where you look among the numbers you find plenty of primes. From a scientific point of view you can do billions of experiments with your mathematical package and always find primes. It can be tested more often and more precisely than any law of physics. You can safely bet the family farm on this and still sleep soundly at night. And yet, I think many of you will agree with me that in this case scientific observation and experimental evidence is a sorry excuse for real knowledge. It may be acceptable for a court of law or everyday life, but it is totally unacceptable given that you can use pure logical reasoning and a few basic axioms for numbers to conclude that this is not just an empirical observation, but a fact built into the structure of whole numbers themselves. It was the genius of the ancient Greeks to develop mathematics not just as an empirical science, but as an axiomatic system for logical deductions.\footnote{The notion of proof may go back much further to the Babylonians.} In Euclid we find in place of the above scientific observation the following deduction.

\bigskip 
\noindent\textbf{Theorem.}  \emph{There are infinitely many prime numbers.} 
\bigskip

\emph{Proof.} Given primes $p_1, p_2, \ldots , p_n$, the number
\begin{equation} N  = (p_1 \cdot p_2 \cdot p_3 \cdots p_n) +	1 \label{3}\end{equation}
must contain a prime factor not among the primes used in its construction. To see this, notice $p_1$ does not
divide $N$ since it leaves a remainder of 1 (or alternatively $N/p_1$ is clearly not an integer). Similarly the other $p_i$'s do not divide $N$. We therefore conclude that any finite list of primes is not complete, and therefore there must be infinitely many primes. As an example, if we start with the primes $2$, $3$, and $29$, we find $N= 2\cdot 3\cdot 29 +1 = 175 = 5^2\cdot 7$ and thus $N$ contains the two new primes $5$ and $7$. 

Now consider the sequence
\begin{equation} 3,5,7,11,13,17,19,29,31,41,43, 71,73,101, 103,\ldots\ . \label{4}\end{equation}
These are also prime numbers, but only the prime numbers that occur as pairs of primes  that are two apart. This is the sequence of twin primes. Once again we might wonder if twin primes continue to occur or eventually die out. Looking for twin primes starting at a billion, we already see one pair in \eqref{2}, and we readily find
\begin{equation}  \begin{split}&1000000007, 1000000009, 1000000409, 1000000411, 1000000931, 1000000933, \\& 1000001447, 1000001449,1000001789, 1000001791, 1000001801, 1000001803 .  \end{split} \label{5}  \end{equation}
There is nothing special about a billion, and you can check that you always find twin primes, although for numbers with hundreds of digits they do not occur nearly as frequently as the prime numbers. This evidence  points strongly to the following conclusion.

\bigskip 
\noindent\textbf{Scientific or Empirical Observation.}  \emph{There are infinitely many twin primes.} 
\bigskip

This is such a natural observation that it is hard to believe that the Greeks did not discover it. Strangely however, the first known published reference to this question was made by  A. de Polignac in 1849, who conjectured that there will be infinitely many prime pairs  with any given even difference. Once again, empirically one can sleep soundly after betting the farm that this observation is true, but unlike for the infinitude of primes, no one has found a string of logical reasoning that demonstrates its truth is built into the structure of the integers. Mathematicians like challenges, and often give names to challenging unsolved problems. 

\bigskip 
\noindent\textbf{Twin Prime Conjecture.}  \emph{There are infinitely many twin primes.} 
\bigskip

You are welcome to try to prove this conjecture and become famous, but be warned that a great deal of effort has already been expended on this problem. The chances that a simple idea such as \eqref{3} will work is very small. Therefore also put some effort into understanding what has been learned about primes in the last two hundred years.

One final word on my argument that empirical evidence provides a solid basis for deciding mathematical questions. One can point to many counterexamples for this, but I think if you examine these you will usually find they occur either because the observations were not sufficient, or because the phenomena had abrupt changes in behavior.  The twin prime conjecture could fail if properties of very large numbers, say with more than a million digits, are vastly different than smaller numbers.
However, since the properties that generate the integers are in play from the start, it is against everything we know to believe that \emph{all} large numbers will behave fundamentally differently than smaller ones. Other famous unsolved problems depend on \emph{every} large number not having an unusual and unlikely property, and in that situation one is on much shakier ground.\footnote{One such problem is the existence of Landau-Siegel zeros, which we will encounter later in this paper.} 

\section{Primes Thin Out }

When people first begin to study prime numbers, they frequently want to find a \emph{formula} for them. Ideally one would plug $n$ into this formula, and the formula would produce the $n$-th prime. While there actually are such formulas (see for example \cite{HR}), they are so complicated to compute with that you are better off using the Sieve of Eratosthenes, a procedure for quickly finding all the primes up to a given size. To find all the primes up to $N$, for example,  you first remove all the multiples of 2 --- the even numbers, from your list of natural numbers up to $N$, then all multiples of 3, then multiples of 5, and so on.  After removing all the multiples of the initial primes $2, 3, 5, 7, \ldots P$,  the numbers left in your list greater than $P$ and less than $P^2$ will be exactly the primes in this range, since the primes in this range will not have been removed, while any composite number with no prime factors $\le P$ must neccesarily be larger than $P^2$. If you therefore pick $P$ to be the largest prime less than $\sqrt{N}$, you will obtain a list of the primes up to $N$. For example, to find all the primes less than a million, you only need to sieve with the 168 primes less than a thousand. Despite this simple procedure for generating the primes, the individual primes appear in a highly irregular pattern,  as a glance at the list in \eqref{2} reveals. This irregularity makes it unrealistic that the primes are obtainable by a simple plug-in formula. 

In view of these considerations, the first step in understanding primes is to think of them as a natural phenomena and look for statistical rather than exact data. Clearly they get rarer as we move towards larger numbers, and there is a simple reason for this. To start with, after 2 all the multiples of 2, the even numbers, can not be primes, and therefore we have eliminated one half of all the natural numbers. 
Next, every multiple of 3, (which could be called the \lq\lq threeven" numbers but have no name in English), can not be primes, and this eliminates a third of the numbers, however half of these multiples of three are also even, and having already been eliminated need to be added back as a correction for this overcounting. Thus the proportion of natural numbers not divisible by 2 and 3 is
\[ 1 -\frac12 -\frac13 + \frac16 = \frac13. \]
You can easily check for yourself that  the proportion of natural numbers not divisible by 2, 3, and 5 is by what is called an inclusion-exclusion process
\begin{equation} 1 -\frac12 -\frac13  - \frac15 + \frac1{2\cdot 3} + \frac1{2\cdot 5}+ \frac1{3\cdot 5}- \frac1{2\cdot 3 \cdot 5 }= \frac4{15} .\label{6}\end{equation}
The sum on the left-hand side is actually
\begin{equation} \left(1-\frac12\right)\left(1-\frac13\right)\left(1-\frac15\right),\label{7}\end{equation}
and now we can see that the above analysis is better understood in terms of simple probability. The probability that a number is not even is $1/2 = (1-1/2)$. The probability that a number is not divisible by 3 is $2/3 = (1-1/3)$, and the probability  that a number is not divisible by 5 is $4/5= (1-1/5)$. But these events are independent of each other, since knowing only that a number has a prime factor of 2 tells you nothing about whether
it has a prime factor of 3 or 5. Hence the event where all these conditions hold is the product of the individual probabilities. Thus the probability that a natural number is not divisible by all the primes up to P is
\begin{equation} \prod_{p\le P}\left(1-\frac{1}{p}\right) = \left(1-\frac12\right)\left(1-\frac13\right)\left(1-\frac15\right)\left(1-\frac17\right) \cdots \left(1-\frac1P\right).\label{8}\end{equation}
Clearly this product is never negative, and  you might guess
\begin{equation}  \lim_{P \to \infty} \prod_{p\le P}\left(1-\frac{1}{p}\right) =0 ,\label{9}\end{equation}
so that the probability a random integer is a prime decreases to zero as the size of the integer gets larger.\footnote{ Numerically one does observe this, although the rate that this proportion goes to zero is surprisingly slow. For example, sieving by the 25 primes less than 100 leaves 12\% of the integers, while sieve out the 1229 primes less than 10,000 still leaves 6\% of the natural numbers. }
This was proved by Legendre.  You might try to prove this yourself, although it is difficult if you haven't see the proof before.  (As a first step, if you take logarithms of both sides of \eqref{9}, the problem can be reduced with a little calculus to showing that
\begin{equation}  \sum_p \frac1p =  \frac{1}{2}+ \frac{1}{3}+\frac{1}{5}+ \frac{1}{7}+\frac{1}{11} + \cdots =\infty;\label{10} \end{equation}
the sum of the reciprocals of the primes diverge. 
You can find a proof in many beginning number theory books, usually in the section on Mertens' theorem. (See for example \cite{CM}, or \cite{HW}.)

There is a standard formulation of the Sieve of Erastothenes involving the M\"obius function $\mu$ defined by
\begin{equation} 
	\mu(n) = 
		\left\{ 
		\begin{array}{ll}
       	1, &\mbox{if $n= 1$}, \\
                  (-1)^r, & \mbox{if $n=p_1 p_2 \cdots p_r$} , \\
                  0, & \mbox{if $n$ has a repeated prime factor}.
		\end{array}
		\right. 
	\label{11}\end{equation}
The product in \eqref{8} can now be expressed as a  sum by 
\begin{equation}\prod_{p\le P}\left(1-\frac{1}{p}\right)= \sum_{d| \mathcal{P}}\frac{\mu(d)}{d}, \quad \mathcal{P} = 2\cdot 3\cdot 5\cdots P,\label{12}\end{equation}
where the condition $d|\mathcal{P}$ means that we sum over all divisors of $\mathcal{P}$. You can see that equation \eqref{6} is an example of this formula with $\mathcal{P}= 2\cdot 3\cdot 5$. We now rewrite \eqref{9} as
\begin{equation} \lim_{P\to \infty} \sum_{d| \mathcal{P}}\frac{\mu(d)}{d} =0. \label{13}\end{equation}
Since this sum will eventually run through all the natural numbers, it would seem clear that
\begin{equation} \sum_{n=1}^\infty \frac{\mu(n)}{n} =0. \label{14}\end{equation}
However, this turns out to be very difficult to deduce\footnote{This should not be too surprising for those of you who have studied  conditionally convergent series in calculus. It is actually the existence of a limit which is hard to prove; it is relatively easy to prove that if there is a limit it must be zero.}, and was only proved in 1899 by Landau, about a hundred years after Legendre's formula. Equation \eqref{14} is  essentially equivalent to the prime number theorem which we will  introduce in the next section.  The failure of sieve methods to prove results like \eqref{14}   led to the dominance of analytic methods in the study of primes for many years. However, many of the recent important advances in the subject have depended on the combination of both sieve and analytic methods.

\section{The Prime Number Theorem}
\begin{figure}
\centerline{\includegraphics{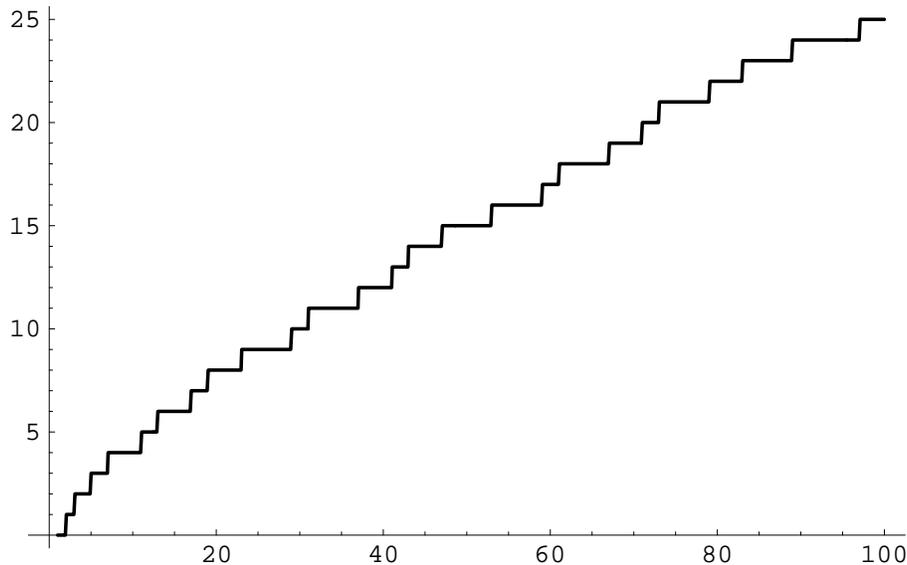}}
\caption{The graph of $\pi(x)$ for $1\le x \le 100$}\label{fig1}
\end{figure}
\begin{figure}
\centerline{\includegraphics{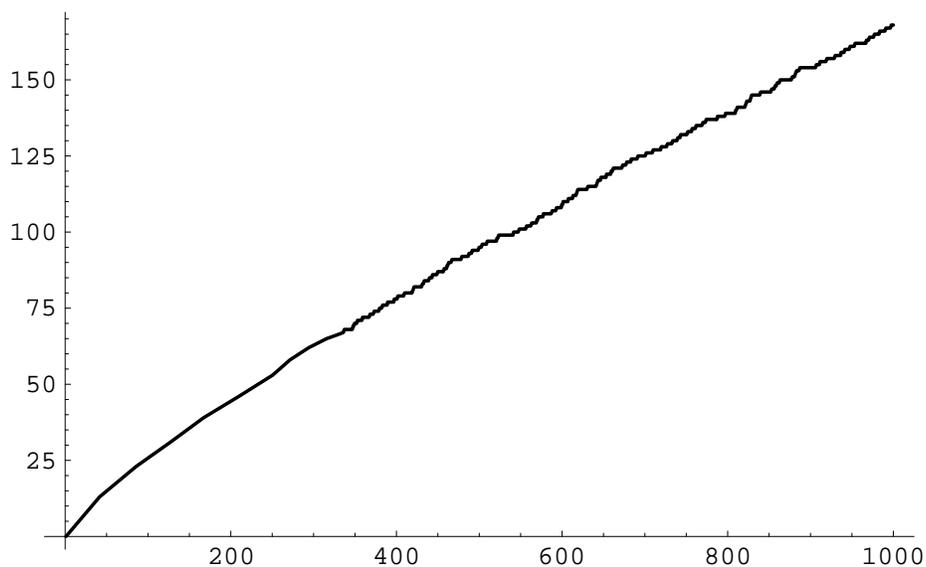}}
\caption{ The graph of $\pi(x)$ for $1\le x \le 1000$}\label{fig2}
\end{figure}

Let us now examine statistically the rate at which the primes thin out. We define $\pi(x)$ to be the number of primes less than or equal to $x$. Thus for example $\pi(5)=3$ since 2, 3, and 5 are the three primes less than or equal to 5. You can see from \eqref{1} that $\pi(100) = 25$. In Figure 1 above is the graph of $\pi(x)$ for $1\le x \le 100$.
Clearly $\pi(x)$ is a step function with jumps at the primes, and therefore when looked at closely is as complicated as the primes. 
But when you move back and view it over a longer range
$\pi(x)$  becomes extremely regular, as you can see in Figures 2 and 3 where $\pi(x)$ is graphed over the range $1\le x \le 1000$, and $1\le x \le 1,000,000$. 
\begin{figure}
\centerline{\includegraphics{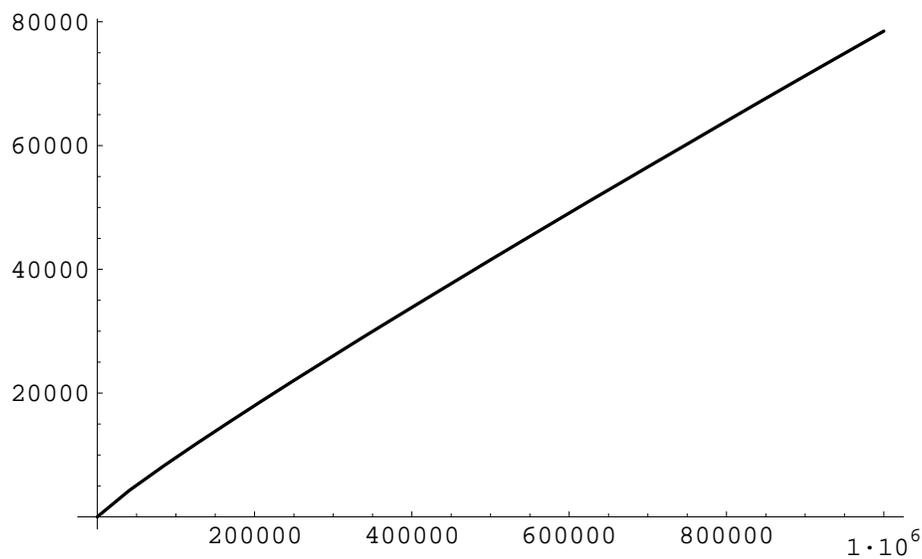}}
\caption{ The graph of $\pi(x)$ for $1\le x \le 1,000,000$}\label{fig3}
\end{figure}

What should be clear from these graphs is that $\pi(x)$ must approach a very simple function as $x\to \infty$. This fact is now called the prime number theorem, and it asserts that
\begin{equation} \pi(x) \sim  \frac{x}{\log x}, \quad \text{as} \ \ x \to \infty, \label{15}\end{equation}
where $\sim$ means the ratio of the two sides approaches 1 as $x\to \infty$. Here $\log$ is the natural logarithm\footnote{The majority of mathematicians use $\log$  instead of $\ln$ to denote the natural logarithm except when they are teaching a calculus class.} to the base $e=2.71828\ldots$, and you might wonder why the primes know about this transcendental number. The answer is because the integers also know about $e$ since 
\[ \sum_{n\le N} \frac1n = 1 + \frac12 + \frac13 +\frac14 +\frac15 +\frac16 + \cdots +\frac1N \sim \log N;\]
the harmonic series grows like the natural logarithm. While this is easily proved by calculus since the series is well approximated by the corresponding integral, the connection with the prime number theorem is much more difficult to establish. 

In Figure 4 we have graphed $\pi(x)$ and $x/\log x$ together; $\pi(x)$ is the upper curve. As you can see, this isn't a very good fit, but it suggests (at least in hindsight) the correct approximation. The prime number theorem says that on average the probability that a number between $1$ and $x$ is a prime is $1/\log x$,  and therefore an individual number $n$ should have  probability $1/\log n$ of being a prime. This density function no longer depends on the large global variable $x$, and to find the total number of primes we should add up the local probabilities or equivalently integrate the density function. Therefore it makes sense to approximate $\pi(x)$ with the logarithmic integral
\begin{equation} \text{li}(x) := \int_2^x \frac{1}{\log t}\, dt \sim \sum_{2\le n\le x}\frac{1}{\log n} .\label{16} \end{equation}
\begin{figure}
\centerline{\includegraphics{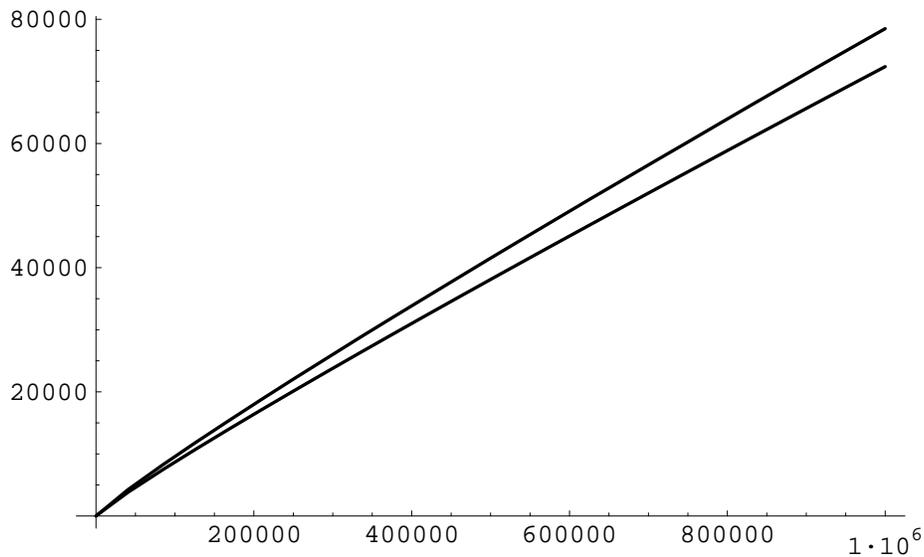}}
\caption{ The graph of $\pi(x)$ (above) and $x/\log x$ (below) for $1\le x \le 1,000,000$}\label{fig4}
\end{figure}
\begin{figure}
\centerline{\includegraphics{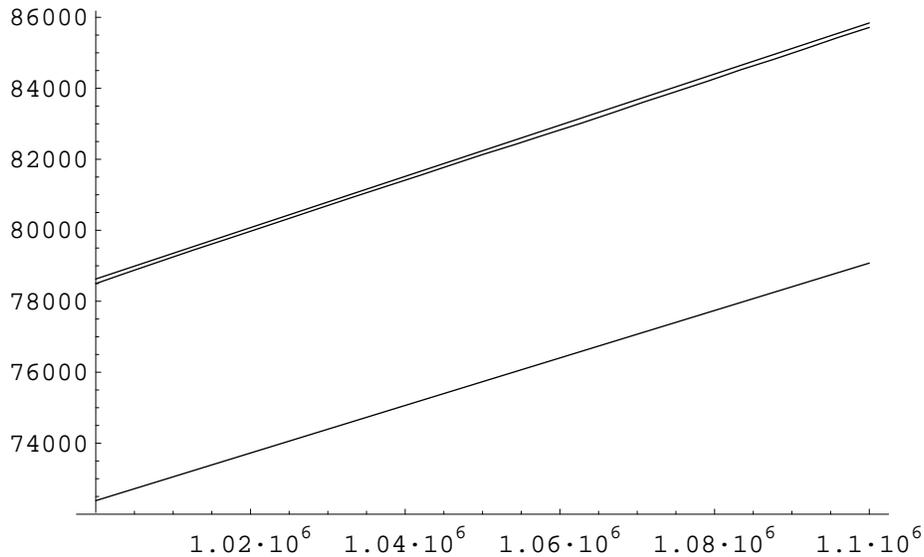}}
\caption{ The graph of $\text{li}(x)$ (top), $\pi(x)$, and $x/\log x$ (bottom) for $1,000,000 \le x \le 1,100,000$}\label{fig5}
\end{figure}
There is no point in comparing the graphs of $\pi(x)$ and $\text{li}(x)$ over the range \mbox{$1\le x\le 1,000,000$} since the two graphs are so close they appear to be the same graph. In Figure 5 we graph  $\pi(x)$, $\text{li}(x)$, and $x/\log x$ in the range 1,000,000 $\le x\le $1,100,000; the top curve is $\text{li}(x)$ and the bottom curve is $x/\log x$. 

This evidence for the prime number theorem was first noticed independently by Gauss and Legendre near the end of the 18th century\footnote{The teenage Gauss found the correct approximation $\text{li}(x)$, while Legendre found $x/\log x$ and speculated incorrectly on higher order approximations.} but its proof was only found in the closing years of the 19th century, when in 1896 Hadamard and de la Vall\'ee Poussin independently proved the prime number theorem. Using integration by parts it is easy to see
\begin{equation} \text{li} (x) = \frac{x}{\log x} + \frac{x}{(\log x)^2}+ \frac{2!x}{(\log x)^3} + \frac{3!x}{(\log x)^4}+ \cdots ,\label{17}\end{equation}
and in 1899 de la Vall\'ee Poussin proved that $\text{li}(x)$ is a better fit to $\pi(x)$ than any finite truncation of the series in \eqref{17}.

It was also noted empirically that $\text{li}(x)$ is always found to be larger than $\pi(x)$, which suggests the conjecture that this will always remain true.  This however turns out to be false, as proved by Littlewood in 1914. This famous result is frequently cited as an example of the danger of using empirical data instead of proofs. However, I think the graphs above make it  rather speculative to guess anything long range for this finer order behavior. 

\section{The Riemann Hypothesis}

The extraordinarily good fit between $\pi(x)$ and $\text{li}(x)$, far better than the first approximation $x/\log x$, has been the  subject of intensive but largely unsuccessful investigation for the last one hundred years. From probability considerations one might expect that the fit should be about (or a little bigger) than the square root of the approximation.  This may be observed empirically for primes, and suggests that the error in the prime number theorem satisfies the bound, for $x$ sufficiently large, 
\begin{equation} \big| \pi(x) - \text{li}(x)\big| < C x^{ \frac12 + \epsilon}, \quad \text{for any} \ \epsilon > 0, \label{18} \end{equation}
and some constant $C$.
This statement turns out to be equivalent to a conjecture Riemann made in 1859 concerning the zeros of the Riemann zeta-function. The Riemann zeta-function is defined by
\begin{equation} \zeta(x +iy) = \sum_{n=1}^\infty \frac{1}{n^{x+iy}}, \quad \text{for} \ \ x>1,\label{19}
\end{equation}
where $i =\sqrt{-1}$. This function can be studied as a function of the complex variable $z=x+iy$. The starting point for the proof of the prime number theorem is the Euler product identity
\begin{equation}\begin{split} \zeta(z) &= \sum_{n=1}^\infty \frac{1}{n^{z}} = \prod_{p}\left(1 + \frac{1}{p^z} + \frac{1}{p^{2z}}+\frac{1}{p^{3z}}+ \cdots \right) \\& =\prod_{p}\left(1-\frac{1}{p^z}\right)^{-1} ,\label{20}\end{split} \end{equation}
which you can verify by seeing how the terms in the first product with primes $2, 3, 5, \ldots$  can be multiplied out to give each natural number exactly once. (The series in the first product is a geometric series which converges for $x>1$ to the expression in the second product.)

The Riemann zeta-function can never be zero if $x>1$. The series and product in \eqref{20} only converge for $x>1$, but one can find alternative expressions which agree with these formulas for $x>1$ but continue to be valid for all $z$ except $z=1$ which is a singularity of the zeta-function.  This extension is the unique differentiable extension, which is called the analytic continuation.  In the vertical strip in the complex plane $0\le x\le 1$ the zeta-function is  equal to zero at many values of $z$; these places are called \lq \lq zeta-zeros " of briefly \lq \lq zeros ". The Riemann Hypothesis is that all these  zeros  occur at complex numbers $z= 1/2 +iy$ on the vertical line with real part equal to 1/2. While it has been verified that the first 10 trillion zeros in this strip above and below the real axis lie on this line, a proof is not in sight. The Clay Institute has offered a million dollar prize for a proof of the Riemann Hypothesis, and also a million dollar prizes for 6 other \lq \lq Millennium " problems. While the Riemann Hypothesis is decisive in determining the distribution of primes, it seems to be of  of little help with regard to twin primes.

\section{The Twin Prime Number Theorem}
What about twin primes? One immediately notices that twin primes thin out faster than primes. Here we have a famous theorem of Brun from 1919 that in contrast to \eqref{10}
\begin{equation} \sum_{\substack{ p\\ p\ \text{ a twin prime}}} \frac1p = \frac13 +\frac15 +\frac17+\frac1{11}+\frac1{13} + \frac1{17}+\frac1{19}+\frac1{29}+\frac1{31}+ \cdots  < \infty ;\label{21}\end{equation}
 the sum of reciprocals of the twin primes converge. 
\footnote{ The sum converges to $1.92016\ldots$. One odd aspect of this type of result is that while we do not know if this sum is over a finite or infinite number of twin primes, we can still compute its value as precisely as available computer resources allow.     It was in computing this constant in 1994 that Nicely discovered a flaw in the Intel Pentium chip, known but not reported by Intel, which created a public furor --- \lq\lq Intel: quality is job 0.999999998", which ended up costing Intel hundreds of millions of dollars in recalled chips.} 

Let $\pi_2(x)$ be the number of twin prime pairs with the smaller prime in the pair less than or equal to $x$. Thus for example $\pi_2(10)=2$ because of the two twin prime pairs $(3,5)$ and $(5,7)$. In Figures 6, 7 and 8 we graph $\pi_2(x)$ just as we did before for $\pi(x)$. While somewhat more irregular, we once again see an asymptotic behavior developing. 
\begin{figure}
\centerline{\includegraphics{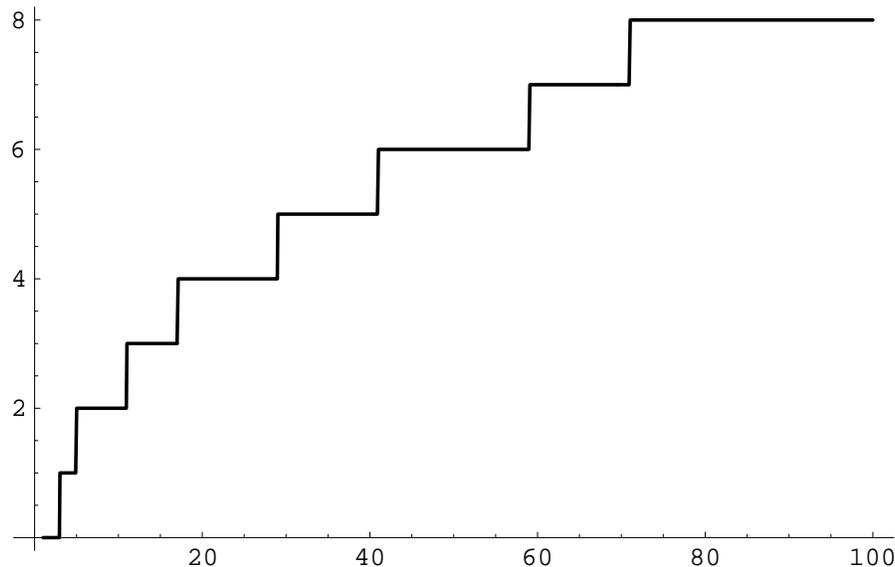}}
\caption{ The graph of $\pi_2(x)$ for $0\le x \le 100$}\label{fig6}
\end{figure}
\begin{figure}
\centerline{\includegraphics{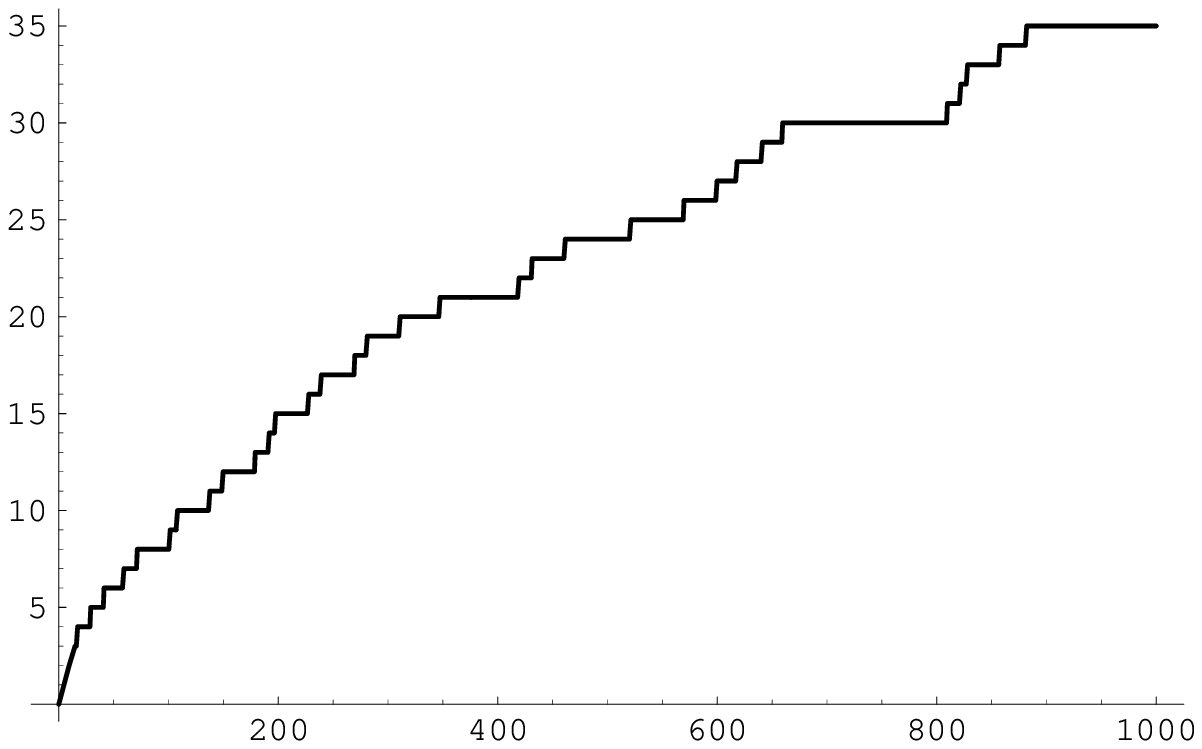}}
\caption{ The graph of $\pi_2(x)$ for $2\le x \le 1000$}\label{fig7}
\end{figure}
\begin{figure}
\centerline{\includegraphics{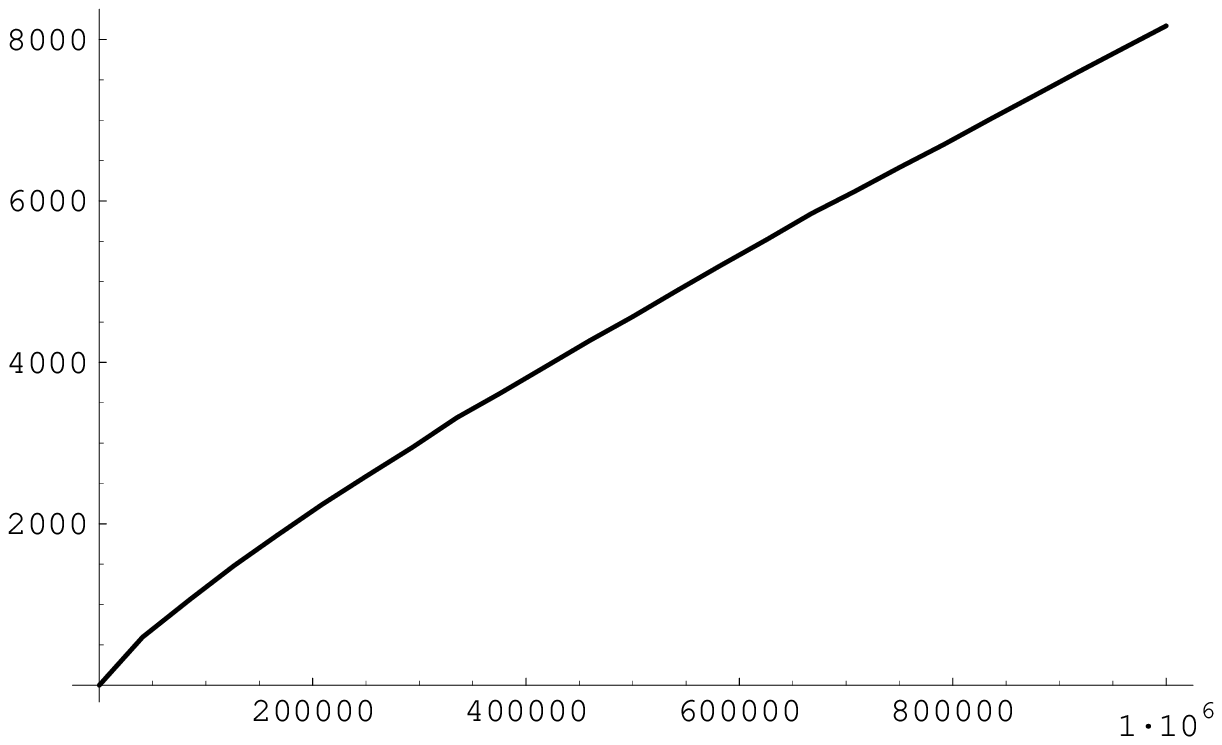}}
\caption{ The graph of $\pi_2(x)$ for $2\le x \le 1,000,000$}\label{fig8}
\end{figure}
\begin{figure}
\centerline{\includegraphics{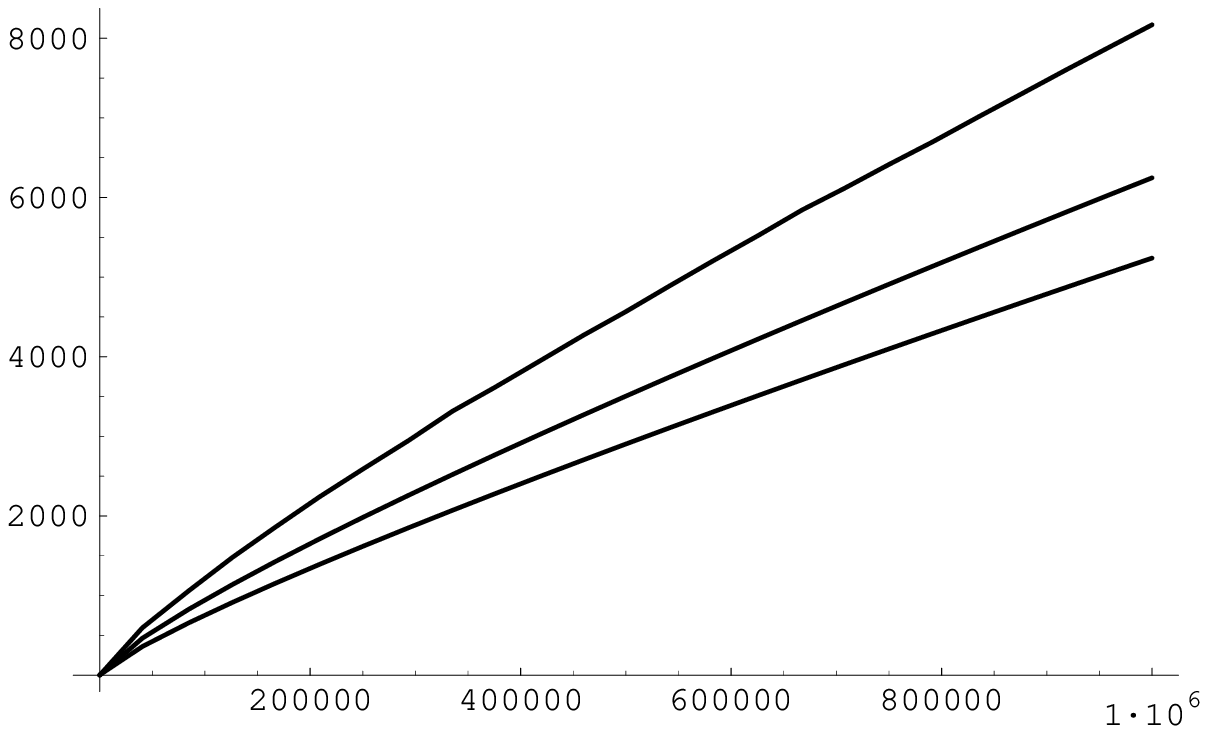}}
\caption{ The graph of $\pi_2(x)$ (top), $\text{li}_2(x)$ (middle), and $x/(\log x)^2$ for $1\le x \le 10^6$}\label{fig9}
\end{figure}

What is the correct asymptotic function?  Returning to probability considerations, what is the probability that $n$ and $n+2$ are both prime? The prime number theorem is consistent with assigning the probability that a random number $n$ is prime to be $1/\log n$; this is called the Cram\'er model, introduced and made use of by H. Cramer in 1935 \cite{Cr}.   The chance that $n+2$ is prime is  then $1/\log(n+2) \sim 1/\log n$. Therefore 
by independence the probability of both being prime is $1/(\log n)^2$.  This suggests that the correct first order approximation for $\pi_2(x)$ should be  $x/(\log x)^2$, or more precisely
\begin{equation}  \text{li}_2(x) := \int_2^x \frac{1}{(\log t)^2}\, dt. \label{22} \end{equation}
We have graphed these in Figure 9, and as you can see we clearly have the wrong answer. 

Although we can not prove anything, there is a heuristically argument which suggests the correct answer.  This argument can be formulated in various ways, but here we follow Soundararajan \cite{So}. One problem with the Cram\'er's model is that it fails to take into account divisibility. Thus, for the primes $p> 2$, 
the probability that $p+1$ is prime is not 
\mbox{$1/\log (p+1)$} as suggested by the Cram\'er model but rather $0$ since $p+1$ is even. 
Further $p+2$ is necessarily odd;  therefore it is twice as likely to be prime as a random number. The conclusion is that $n$ and $n+2$ being primes are not independent events. Let us now correct for this lack of independence.
First, for large twin primes, we need both $n$ and $n+2$ to not be divisible by $2,3,5,7, 11, \cdots$.
The chance that two random numbers are both odd is $(1/2)(1/2) = 1/4$, but, since $n$ being odd forces $n+2$ to be odd, the chance that $n$ and $n+2$ are both odd is $1/2$, and thus twice as large as random.
The chance that two random numbers are both not divisible by 3 is $(2/3)(2/3) = 4/9$,
but the chance that $n$ and $n+2$ are not both divisible by 3 is $1/3$ since this occurs if and only if $n$ is congruent to 2 modulo 3.

In general, the probability that two random numbers are not divisible by $p>2$ is $(1-1/p)^2$,
while the probability that both $n$ and $n+2$ are not divisible by $p$ is the slightly smaller $1 - 2/p$ 
since $n$ must miss the two residue classes $0$ and $-2$ modulo $p$.
Therefore, the correction factor to the Cram\'er model for lack of independence is $2$ if $p=2$, and for $p\ge 3$ is
\[\begin{split} \left(1-\frac2p\right)\left(1-\frac1p\right)^{-2}& = \left(\left(1-\frac1p\right)^2 - \frac1{p^2}\right) \left(1-\frac1p\right)^{-2}\\& =\left(1 -\frac{1}{(p-1)^2}\right).\end{split}\]
We conclude that the correct approximation for twin primes should be
\begin{equation} \pi_2(x) \sim 2\prod_{p>2}\left(1 -\frac{1}{(p-1)^2}\right) \frac{x}{(\log x)^2}. \label{23} \end{equation}
Given our experience from the prime number theorem, we formulate this as the following conjecture.

\bigskip 
\noindent\textbf{The Twin Prime Number Theorem Conjecture.}  \emph{We have} 
\begin{equation} \pi_2(x) \sim 2\prod_{p>2}\left(1 -\frac{1}{(p-1)^2}\right) \mathrm{li}_{\, 2}(x). \label{24}\end{equation}
\bigskip

The twin prime constant has been computed to many digits, and it is known that 
\begin{equation} 2\prod_{p>2}\left(1 -\frac{1}{(p-1)^2}\right) = 1.32032362\ldots .\label{25}\end{equation} 
In Figure 10 we see that this conjecture provides what is surely the correct approximation to $\pi_2(x)$. One might even conjecture that this approximation holds with a square root error. The conjecture that the distribution of twin primes satisfy a Riemann Hypothesis type error term is well supported empirically, but I think this might be a problem that survives the current millennium.
\begin{figure}
\centerline{\includegraphics{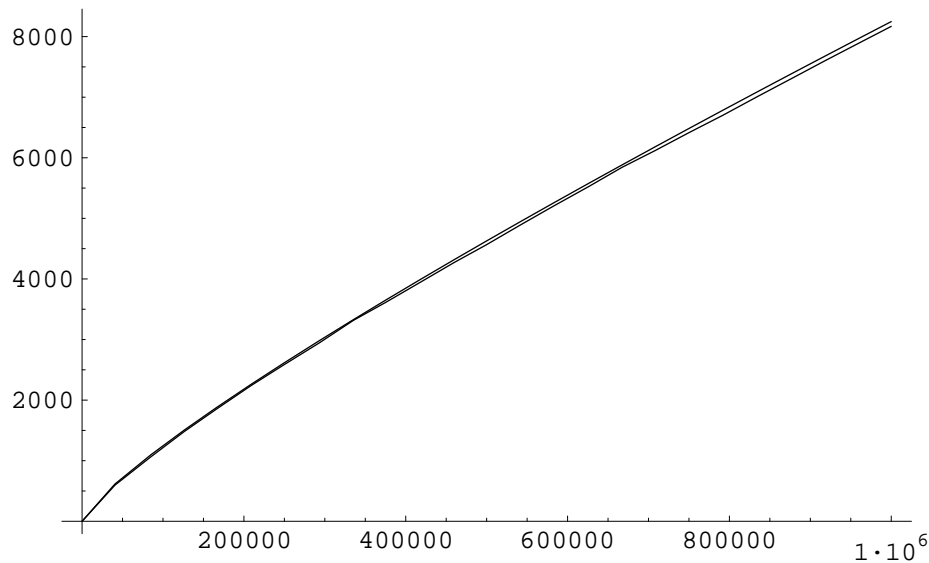}}
\caption{ The graph of $\pi_2(x)$ (bottom) and $ 1.32032362 \ \text{li}_2(x)$ (top) for $1\le x \le 1,000,000$}\label{fig10}
\end{figure}
 
\section{ Some recent progress towards the twin prime conjecture.}

One appealing aspect of number theory is that it is hard to predict which problems can be solved at our present state of knowledge, and which are currently beyond hope of solution. Fermat's Last Theorem, for example, looked hopelessly hard until a totally new approach was discovered, and even that approach had extremely difficult obstacles. And yet in 10 years the proof was complete.

In the case of twin primes, there have been several outstanding advances which to outward appearances could seem just as formitable as the twin prime conjecture.  To mention just two of these results, it has been proved that  for large enough $x$,
\begin{equation} \pi_2(x) \le 8\prod_{p>2}\left(1 -\frac{1}{(p-1)^2}\right) \frac{x}{\log^2 x} .\label{26}\end{equation}
Comparing this with \eqref{23} or \eqref{24} we see that there can be at most 4 times as many twin primes as conjectured.
Secondly, J. Chen proved that there are infinitely many primes $p$ where $p+2$ is either a prime or a product of two primes, see for example \cite{HR}. 

Neither of these results address the question of finding prime numbers which are close together. This problem has a long history, and I would like to conclude this paper by mentioning some new work of Pintz, Y{\i}ld{\i}r{\i}m and me on this topic.  The question we were initially investigating was not directed at the twin prime conjecture, but rather the question of finding smaller than average gaps between primes.  By the prime number theorem the average distance between two consecutive primes in the interval $[0,x]$ is
\begin{equation} \frac{ \text{ length of } [0,x]}{ \text{number of primes in }[0,x]}= \frac{x}{\pi(x)} \sim \frac{  x}{\frac{x}{\log x}} \sim \log x .\label{27}\end{equation}
Thus for example around $x=10^6$ primes are on average about $\log (10^6) = 6\log(10) =13.81\ldots \ $ apart and this spacing doubles to $27.63\ldots$ at $x=10^{12}$. We can now ask if there are always going to be primes substantially closer than this average as $x$ gets larger and larger. To examine this question, consider the sequence, where $p_n$ denotes the $n$-th prime,
\[  \left\{ \frac{p_{n+1}-p_n}{\log p_n}\right\}_{n=1}^\infty ;
\]
we expect that these values will infinitely often be small. Mathematically we measure this by looking for the smallest limit point of the sequence, i.e. the  \lq \lq lim inf".\footnote{Erd\"os proved that this sequence has many limit points in the sense that the set of limit points has positive Lebesgue measure. Unfortunately the proof does not tell us the value of any of these limit points.}
Thus we define
\begin{equation}\Delta := \liminf_{n\to \infty} \left( \frac{p_{n+1}-p_n}{\log p_n}\right)
\label{28}\end{equation}
and hope to prove $\Delta$ is small. Remarkably, 80 years of work had only found that $\Delta\le .248\ldots$, a result of Maier \cite{Ma} from 1988 which utilized all the previous methods applied to the problem. Then last year we finally were able to prove the result suggested by the twin prime conjecture.

\noindent \textbf{Theorem.} (Goldston-Pintz-Y{\i}ld{\i}r{\i}m) We have $\Delta =0$.

\noindent Our method thus produces primes very close together in a statistical sense, but what came as a great surprise to us is that if you assume an unproved conjecture concerning primes in arithmetic progressions then the method actually produces primes that are a bounded distance apart. That one can prove such a strong result using this information runs counter to all previous expectations, and for several weeks this convinced us that there must be a mistake in our proof. 

To describe this result we begin with a simple example. If 
you divide the natural numbers up modulo 3 you get three residue classes or arithmetic progressions:

\[n\equiv 0 (\text {mod }3):\   \ 3, 6, 9,12,15 \ldots ,\]
\[n\equiv 1 (\text {mod }3):\   1,4,7,10,13, \ldots , \] 
\[ n\equiv 2 (\text {mod }3):\   2,5,8,11,14, \ldots .\]

\noindent Clearly the only prime in the progression $0 (\text{mod} \, 3)$ is 3, but we expect that the primes should be equally distributed in the other two progressions, and hence, letting
$\pi(x;q,a)$ denotes the number of primes $\le x$ in the progression  $a (\text{mod} \, q)$,
\[ \pi(x;3,1) \sim \pi(x;3,2) \sim \frac{1}{2}\ \pi(x) \sim \frac{1}{2} \frac{x}{\log x}. \]

The generalization of this result is called the prime number theorem for arithmetic progressions. We need to avoid progressions like $0 (\text{mod}\, 3) $ where each term is a multiple of some integer $\ge 2$;  of the $q$ progressions modulo $q$ the number of progressions which are not multiples of some number is $\phi(q)$, the Euler phi function, defined by
\begin{equation} \phi(q) := \# \big\{a: 1\le a \le q  \ \text{and} \  (a,q)=1 \big\}, \label{29} \end{equation}
where here the notation $(a,q)$ is the gcd of $a$ and $q$. 
The prime number theorem for arithmetic progressions states that for $(a,q)=1$,
\begin{equation} \pi(x;q,a) \sim \frac{1}{\phi(q)}{\rm li}(x) . \label{30}\end{equation}
For applications we need $q=q(x) \to \infty$ as $x\to \infty$, but unfortunately the best result known allows $q$ to only grow at the very slow rate 
$q\le (\log x)^A$, for any $A$.\footnote{The extension of this result to larger $q$ requires the solution of two famous problems involving Dirichlet $L$-functions, which are a class of functions that include the Riemann zeta-function introduced earlier. The first problem is to show that there are no zeros on the real axis, called Landau-Siegel zeros, and the second problem is to prove the Riemann Hypothesis for Dirichlet L-functions.}  
However, in applications it is often enough to know that on average over many progressions the error here is small, and for this one can take $q$ much larger. The main result of this type was proved in 1965 independently by Bombieri and Vinogradov, and states that
for any $A>1$ we have
\begin{equation}
\sum_{q\le Q}  \max_{\substack{a\\ (a,q)=1}}  \left| \pi(x;q,a)-   \frac{{\rm li}(x)}{\phi(q)}  \right| \le C\frac{x}{(\log x)^A} 
\label{31} \end{equation}
for $Q=x^{1/2}/(\log x)^{B}$, where $B$ and $C$ are constants which depends on the given $A$.

The largest power of $x$ which we can take $Q$ to be in the above result is called the \emph{level of distribution} of primes in arithmetic progressions.  Thus the Bombieri-Vinogradov theorem says the primes have level of distribution $\frac12$. More precisely, we define the level of distribution to be $\vartheta$ if \eqref{31} holds for any $\epsilon >0$ and $Q=x^{\vartheta -\epsilon}$. We expect and find numerically that the primes actually have level of distribution $\vartheta=1$; this was first conjectured by Elliott and Halberstam. 
What Pintz, Y{\i}ld{\i}r{\i}m and I proved is that if the primes have a level of distribution equal to any number larger than $\frac12$, then there must be  infinitely often primes a bounded distance apart. In the case of level of distribution 1, we proved the following result.

\medskip
\noindent \textbf{Theorem.} \  If the Elliott-Halberstam Conjecture is true (actually if $\vartheta \geq .971$), then
\begin{equation}  p_{n+1}-p_n \le 16 \quad \text{for infinitely many $n$}.\label{32}\end{equation}

The proof of these results is not that difficult compared to other results in the field, but I can only describe the main ideas here.
In the first place, we need  a generalization of twin primes where we 
consider the tuple or vector
\begin{equation}(n+h_1, n+h_2, \ldots , n+h_k)\label{33}\end{equation} with the shifts $h_i$  given by
$\mathcal{H}= \{h_1,h_2,\ldots , h_k\}.$
Letting $n=1,2,3,\ldots \ $ we ask how often all of the components of the tuple are simultaneously prime for $n\le x$, and denote this number by $\pi(x; \mathcal{H})$. Thus for example  twin primes correspond to $\mathcal{H}=\{0,2\}$ with the tuple $(n,n+2)$. On the other hand 
the tuples $(n,n+1)$ is only made up of primes when $n=2$ since at least one of these numbers is even. Similarly $(n,n+2,n+4)$  is only made up of primes when $n=3$ since at least one of these numbers is divisible by 3. Tuples which do not always have a component divisible by some integer $\ge 2$ are called \emph{admissible}, and for these we expect that infinitely often all their components will simultaneously be primes. This is called the Hardy-Littlewood prime tuple conjecture.
Hardy and Littlewood also made a more precise conjecture. Let $\nu_p(\mathcal{H})$ denote the number of distinct
residue classes$\pmod p$ the numbers $h\in \mathcal{H}$ fall into. Just as for the twin prime constant, we can correct for the lack of independence\footnote{This is done in Soundararajan's paper \cite{So}.} and obtain an expected proportion, called the singular series
\begin{equation}
\mathfrak{S}(\mathcal{H}) = \prod_p\left(
1-\frac{1}{p}\right)^{-k}\left(1 -
\frac{\nu_p(\mathcal{H})}{p}\right).
\label{34}\end{equation}
If $\mathfrak{S}(\mathcal{H})\neq 0$ then $\mathcal{H}$ is admissible.
Thus $\mathcal{H}$ is admissible if and only
if $\nu_p(\mathcal{H})<p$ for all $p$.

\noindent \textbf{Prime Tuple Conjecture} If $\mathcal{H}$ is admissible, then
\begin{equation} \pi(x,\mathcal{H}) \sim \mathfrak{S}(\mathcal{H}) \text{li}_k(x), \ \
\text{where} \ \
\text{li}_k(x) = \int_2^x \frac{1}{(\log t)^k}\, dt .\label{35}\end{equation}

The first idea in our method is to try to replace this counting function by an approximation for which we can prove asymptotic formulas corresponding to the prime tuple conjecture. For this we introduce the 
 von Mangoldt function $\Lambda(n)$, defined to be $\log p$ if $n=p^m$ and zero otherwise. This function tells us whether $n$ is a prime or prime power, but since the number of powers is very small they can be removed from consideration at a later stage.\footnote{You can see easily there are $\le \sqrt{x}$ squares that are $\le x$, and $\le \log_2 x$ different powers $\le x$.}  You might try to prove in a few lines that the prime number theorem in the form \eqref{15} implies and is easily obtained from the formula
\begin{equation}\sum_{n\le x} \Lambda(n) \sim x .\label{36}\end{equation}
In a first course in number theory we prove
 the elementary formula, which you might also try to prove yourself,
\begin{equation}\Lambda(n) = \sum_{d|n} \mu(d) \log \frac{n}{d}.
\label{37}\end{equation}
This formula has little utility in applications because it has too many terms. To understand this last statement, you can  try to prove \eqref{36} by substituting in \eqref{37} and seeing what goes wrong. One can, however, by direct substitution and using the prime number theorem obtain formulas like \eqref{36} for the truncated smoothed approximation 
\begin{equation} \Lambda_R(n) = \sum_{\substack{d|n\\ d\le R}
}\mu(d) \log \frac{R}{d} \label{38} \end{equation}
if $R$ is kept somewhat smaller than $N$.
This approximation may seem ad hoc, but it arises naturally.  
Hardy and Littlewood also formulated the prime tuple conjecture in terms of the von Mangoldt function by defining
\begin{equation} \Lambda(n;\mathcal{H}) := \Lambda(n+h_1)\Lambda(n+h_2)\cdots \Lambda(n+h_k),\label{39}\end{equation}
(so this will be zero if any of the $\Lambda(n+h_i)$ is zero), and then equivalently to \eqref{35} conjectured
\begin{equation} \sum_{n\le N}\Lambda(n;\mathcal{H}) \sim \mathfrak{S}(\mathcal{H})N .\label{40}\end{equation}
In view of \eqref{38} and \eqref{39}, it is natural to approximate $\Lambda(n;\mathcal{H})$ by
\begin{equation}\Lambda_R(n;\mathcal{H}):= \Lambda_R(n+h_1)\Lambda_R(n+h_2)\cdots
\Lambda_R(n+h_k).\label{41} \end{equation}
Until 2004, this was the only useful approximation we knew of for this problem.

We now try to detect primes with these approximations. For this, we need an approximation which is never negative, but the approximation $\Lambda_R(n)$ and also $\Lambda_R(n,\mathcal{H})$ is frequently negative. (For example you can check that $\Lambda_5(30) = -2\log 5$.)  Therefore we need to first square the approximation to obtain a non-negative approximation. The key formulas we need to compute for our method are\begin{equation} \sum_{n\le N}\Lambda_R(n;\mathcal{H})^2 \quad
\text{and}
\quad \sum_{n\le N}\Lambda(n+h_0)\Lambda_R(n;\mathcal{H})^2.\label{42}\end{equation}
While these are complicated to evaluate, the analysis needed is at the level of the prime number theorem with remainder, and therefore classical. For the second sum, the single factor of $\Lambda(n+h_0)$ really is detecting primes (and prime powers) in the tuple, since we find that if $h_0\in \mathcal{H}$ we get a result larger by a factor of $\log R$ than we get when $h_0 \not \in \mathcal{H}$. The second sum is evaluated by summing $\Lambda(n+h_0)$  through arithmetic progressions modulo products of divisors from $\Lambda_R(n;\mathcal{H})^2$, and this is where the level of distribution information is used. The result of this analysis is that for $R = N^{1/4k-\epsilon}$ we obtain asymptotic formulas for both sums in \eqref{42}. Using these formulas, we can evaluate asymptotically, with $r\ge 1$, 
\begin{equation}\mathcal{S} := \sum_{n= N+1}^{2N}\left(\sum_{i=1}^k \Lambda(n+h_i) - r\log (3N) \right) \Lambda_R(n,\mathcal{H})^2 \label{43}\end{equation}
Here, if $1\le h_i<N$, we have $\Lambda(n+h_i)\le \log(n+h_i) <\log 3N$. Thus if we find $\mathcal{S}>0$ then there must be an $n$ for which at least $r+1$ of the $\Lambda(n+h_i)\neq 0$, and (after removing  prime powers) there are $r+1$ primes in the tuple $\mathcal{H}$. 

Unfortunately it turns out that $\mathcal{S}<0$ for the approximation in \eqref{41} even when $r=1$, and we fail to prove there are even two primes in a tuple. However we are able to recover something by this method, if we switch to the more modest goal of finding two primes close together. For this, we now try to find as many primes as possible in the interval $(n,n+h]$, which we detect by using \emph{all} the possible $k$-tuples that can be formed in that interval. Thus we now consider 
\begin{equation} \mathcal{S}' := \sum_{n= N+1}^{2N}\left(\sum_{1\le h_i\le h} \Lambda(N+h_i) - r\log (3N) \right)\sum_{\substack{1\le h_1 ,h_2,\ldots h_k\le h \\ \text{distinct}}} \Lambda_R(n,\mathcal{H})^2 \label{44}\end{equation}
With $r=1$ we find this is positive if $h> \frac34 \log N$, and therefore we conclude $\Delta \le 3/4$.
One can now improve on this analysis by using approximations not just for $k$-tuples but a linear combination of all the approximations for $2$-tuples, $3$-tuples, and so on up to $k$-tuples. This leads to an optimization problem which when solved results in $\Delta \le \frac14$. 

The next step is to search for better approximations than \eqref{41}. The main disadvantage of this approximation is that it is formed from many short divisor sums multiplied together; in \eqref{43} there are $2k$ of them, which is what forces the very short truncation length $R=N^{1/4k-\epsilon}$. This reduces the quality of the approximation. It is natural to think that if we could approximate the number of primes in a tuple by a single divisor sum then we could take the approximation much longer and obtain a better result. As it turns out, sieve methods are based exactly on this idea. 
Instead of the tuple $(n+h_1, n+h_2, \ldots , n+h_k)$ consider the polynomial
\begin{equation}\mathcal{P}(n,\mathcal{H}) = (n+h_1)(n+h_2)\ldots (n+h_k) \label{45} \end{equation}
If our tuple is a prime tuple then $\mathcal{P}$ has $k$ prime factors, and conversely. 
The generalized von Mangoldt function
\begin{equation}\Lambda_k(n) = \sum_{d|n} \mu(d) (\log \frac{n}{d})^k \label{46}\end{equation}
is the arithmetic function commonly used to detect whether numbers have $\le k$ distinct prime factors. It can be proved that $\Lambda_k(n)$ is zero if $n$ has more than $k$ distinct prime factors, but is non-zero if $n$ has $\le k$ prime factors. Therefore
$\Lambda_k(\mathcal{P}(n,\mathcal{H}))$
will be non-zero if the tuple associated with $\mathcal{H}$ is a prime tuple. In view of \eqref{38} it is clear we should approximate this with
\begin{equation}
\Lambda_R(n; \mathcal{H})
= \frac{1}{k!}\sum_{\substack{ d|P_{\mathcal{H}}(n)\\ d\le R}}
 \mu(d) \left(\log \frac{R}{d}\right)^{k } .
\label{47} \end{equation}
(Here the factor $\frac{1}{k!}$ is  a natural normalization.) Notice now that we can approximate a prime tuple with a single divisor sum. This is a big step forward, but when this approximation is used in the previous analysis one still does not find primes in tuples, and one ends up obtaining $\Delta \le .1339\ldots = 1- \sqrt{3}/2$. This was disappointing to Y{\i}ld{\i}r{\i}m and me in 2004, but we should have taken to heart the advice: \lq \lq Never give up!" It turns out only one more idea is needed to break through the barrier, and this was discovered by Pintz. The idea is so simple that when I tell it to you I'm sure you will not believe something like this is what mathematicians get paid for, but I think it usually is the case that good mathematics comes down to common sense reasoning.

The idea is that we have been trying to do too much, when much less would still be much more than we need. No one has ever proved there are infinitely many prime tuples, so as a start we only need to try to find tuples with \emph{SOME} primes in them. For example, if we have a 1000-tuple which has a total of 1500 prime factors in its components, then there must still be at least 500 prime components. Therefore, to detect some primes in tuples, you only need to show that $P_{\mathcal{H}}(n)$ has less than $k+\ell$ prime factors, for some $\ell <k$.   Thus we should consider approximations with $\Lambda_{k+\ell}$, where now we have a new variable $\ell$ to make use of.  Hence, we define our new approximation
\begin{equation}
\Lambda_R(n; \mathcal{H}, \ell)
= \frac{1}{(k + \ell)!}\sum_{\substack{ d|P_{\mathcal{H}}(n)\\ d\le R}}
 \mu(d) \left(\log \frac{R}{d}\right)^{k + \ell} .
\label{48} \end{equation}
This finally succeeds in proving $\Delta =0$ using \eqref{44} with  $k$ large enough and an appropriate choice of $\ell$. And this also succeeds in proving $\mathcal{S}>0$ when $r=1$ in \eqref{43} if the level of distribution $\vartheta >1/2$. If $\vartheta=1$ then we find with $r=1$ that $\mathcal{S}>0$ if $k=7$ and $\ell=1$, so that every admissible 7-tuple has two primes in it infinitely often under this assumption. While these results are a great advance forward, they are still only the camel's nose in the tent, since if $r\ge 2$ we  fail to show $\mathcal{S}>0$ even if the level of distribution is $\vartheta=1$. 

\section{A Curious History}

The history behind the development of the method I have just described is much more convoluted than you might guess. The detection method in \eqref{43} is nearly identical to a method Selberg introduced in 1950 for proving $n$ and $n+2$ will together have 5 or fewer prime factors infinitely often (see \cite{S1992}). This method was generalized by Heath-Brown \cite{HB} in 1997 to $k$-tuples, and this work contains both the detection method and the approximation \eqref{48} in the case $\ell=1$. However the approach was never directly applied  to primes and was never viewed as a possible lower bound method. From 1999 until 2003 Yildirim and I were working on special cases of the approximation \eqref{41} and the formula \eqref{42}. However we saw the problem in terms of probability and approximation of moments, and never considered expressions like $\mathcal{S}$. In 2003 Yildirim and I thought we had proved Theorem 1 with a new approximation somewhat similar to \eqref{47} but more complicated and based partly on guesswork.  We had no idea at the time of the relevance of the generalized von Mangoldt function. In examining our proof, Granville and Soundararajan simplified the original moment method into the form of \eqref{43} and \eqref{44}, and then found a fatal mistake: our approximation did not actually have asymptotic formulas in \eqref{42}. Returning to the old approximation \eqref{41}, we were able to use this new detection method and complete the proof that $\Delta \le \frac14$, but until mid-2004 we had no evidence of any improved approximations, and were ready to believe they did not exist. Also in early 2004 Green and Tao were looking for a special type of sieve bound for primes and numbers with a few prime factors in tuples, but could not find anything appropriate in the literature. Granville brought to their attention a manuscript of our work, and they were able to use the asymptotic formula for the first sum in \eqref{42} in their celebrated work on finding arbitrarily long strings of primes in arithmetic progressions.\footnote{Up until 2004 it had only been proved that there are arithmetic progression of three primes infinitely often. }  The feature they needed in this formula is that each component of the tuple can  be shifted individually in the formulas while the other factors remain unchanged. In sieve methods one does not use expressions like \eqref{41}, but rather expressions like \eqref{47} and \eqref{48}, where changing one component of the tuple changes the entire polynomial $\mathcal{P}_{\mathcal{H}}(n)$. Since we were using a moment method, the approximation \eqref{41} naturally  arises from multiplying out moments and getting products of approximations.  However for small gaps between primes, it is the sieve approximations which are more effective, but we did not realize this until mid-2004.  In retrospect, once one has the Granville and Soundararajan formulation of the problem, it is a small step to move to the sieve type approximations. Ironically, the approximation \eqref{41} has now been discarded in our final work, but perhaps it might still play some future role in the study of twin primes. 

What has been left out of the above account of our work is the frequent contact with other mathematicians who freely provided their ideas and suggestions. Many times these were decisive for getting back on track and moving the work forward. I think in this internet age of quick and easy contact we can take advantage of the experience of the worldwide community of mathematicians in our field, always being careful however to not waste someone's time. 

Finally, while there is no evidence of this left in the final work, at many stages a mathematical package (in my case Mathematica) was the only tool we had for testing ideas and experimenting with guesses. Having not grown up with them it is always an effort for me to use these programs, but there were many stages of the work when I would have quit if not for the information they provided.

\section{Some References for Further Study}

For background information one can in the first place study any beginning number theory book. At a somewhat more difficult level, most elementary results mentioned in this paper can be found in Hardy and Wright \cite{HW}. For the Riemann zeta-function, the classic reference is Titchmarsh \cite{T}. For sieve methods, the recent book of Cojocaru and Murty \cite{CM} makes for interesting reading, while the classical reference is Halberstam and Richert \cite{HR}.  One should also study the long article on sieves by Selberg \cite{S1992}. My favorite research papers related to this area are the 1923 paper of Hardy and Littlewood \cite{HL3}, the 1965 paper of Bombieri and Davenport \cite{BD}, and Selberg's famous 4 page 1947 paper introducing the Selberg sieve. One can find the recent work of Goldston-Pintz-Y{\i}ld{\i}r{\i}m in \cite{GPY}, a short proof in \cite{GMPY}, and a very accessible exposition in \cite{So}. 

\section{Acknowledgement}

I would like to thank Dave and Marilyn Blockus, Brian Conrey, Dashiell Fryer, J\'anos Pintz, Tatianna Shubin, and  Yal\c{c}in Y{\i}ld{\i}r{\i}m for many helpful comments concerning this exposition.

\end{document}